\documentclass[11pt]{article}
\usepackage{mathrsfs}
\usepackage{graphicx}
\usepackage{latexsym,amsmath,amssymb,amsfonts}
\usepackage{color,colordvi,amscd,indentfirst,graphics}
\allowdisplaybreaks

\newtheorem{theorem}{Theorem}[section]

\newtheorem{definition}[theorem]{Definition}
\newtheorem{conjecture}[theorem]{Conjecture}
\newtheorem{claim}[theorem]{Claim}
\newtheorem{corollary}[theorem]{Corollary}

\def\qed{\hfill \rule{4pt}{7pt}}
\def\pf{\noindent {\it Proof.} }

\begin{document}

\title{Spanning trees of $K_{1,4}$-free graphs with a bounded number of leaves and branch vertices }

\author{Pham Hoang Ha\footnote{E-mail address: ha.ph@hnue.edu.vn }\\
	Department of Mathematics\\
	Hanoi National University of Education\\
	136 XuanThuy Str., Hanoi, Vietnam}%

\date{}
\maketitle{}
\bigskip

\begin{abstract}
Let $T$ be a tree. A vertex of degree one is a \emph{leaf} of $T$
and a vertex of degree at least three is a \emph{branch vertex} of
$T$. A graph is said to be \emph{$K_{1,4}$-free} if it does
not contain $K_{1,4}$ as an induced subgraph. In this paper, we study the spanning trees with a bounded number of leaves and branch vertices of $K_ {1,4}$-free graphs. Applying the main results, we also give some improvements of previous results on the spanning tree with few branch vertices for the case of $K_{1,4}$-free graphs.
\end{abstract}



\newpage

\section{Introduction}
In this paper, we only consider finite graphs without loops or multiple edges. 
Let $G$ be a
graph with vertex set $V(G)$ and edge set $E(G)$. For any vertex
$v\in V(G)$, we use $N_G(v)$ and $d_G(v)$  to denote the set of neighbors of $v$ and the degree of $v$ in $G$, respectively. 
We define $G-uv$ to be the
graph obtained from $G$ by deleting the edge $uv\in E(G)$, and
$G+uv$ to be the graph obtained from $G$ by adding an edge $uv$
between two non-adjacent vertices $u$ and $v$ of $G$.
For any $X\subseteq V(G)$, we denote by $|X|$ the cardinality of $X$. 
Sometime, we use $|G|$ to denote $|V(G)|$. We define
$N_G(X)=\bigcup\limits_{x\in X}N_G(x)$ and $\deg_G(X)=\sum\limits_{x\in
	X}\deg_G(x)$. The subgraph of
$G$ induced by $X$ is denoted by $G[X]$.  

A subset $X\subseteq V(G)$ is called an \emph{independent set} of
$G$ if no two vertices of $X$ are adjacent in $G$.  
The maximum size of an independent set in $G$ is denoted by $\alpha(G)$. 
For each positive integer $p$, we define 
$$\sigma_p(G)=\left\{\begin{array}{ll}
+\infty, & \;\mbox{ if } \alpha(G) < p,\\
\min\{\sum\limits_{i=1}^p
d_G(v_i)\;|\;\{v_1,\ldots,v_p\}\, \text{is an independent set in } G\},
& \;\mbox{ if } \alpha(G) \geq p.
\end{array}\right.$$

Let $T$ be a tree. A vertex of degree one is a \emph{leaf} of $T$
and a vertex of degree at least three is a \emph{branch vertex} of
$T$. The set of leaves of $T$
is denoted by $L(T)$ and the set of branch vertices of $T$ is denoted by $B(T)$.

There are several sufficient conditions on the independence
number and the degree sum for  a
graph $G$ to have  a spanning tree with a bounded number of leaves
or branch vertices. Win~\cite{Wi} obtained the following theorem, which confirms a conjecture of Las Vergnas~\cite{LV71}. Beside that, recently, the author \cite{H} also gave an improvement of Win by giving an independence number condition for a graph having a spanning tree which covers a certain subset of $V(G)$ and has at most $l$ leaves. 
\begin{theorem}[{\cite[Win]{Wi}, \cite[Ha]{H}}]\label{t1}
	Let $m\geq 1$ and $l\geq 2$ be integers and let $G$ be a $m$-connected graph.  If
	$\alpha(G)\leq m+l-1$, then $G$ has a spanning tree with at most $l$
	leaves.
\end{theorem}
 As a corollary of Theorem \ref{t1}, we have a sharp result (as a note in \cite{H}) for a connected graph to have a bounded number of branch vertices.
\begin{corollary}
	Let $m\geq 1$ and $k\geq 0$ be two integers and let $G$ be a $m$-connected graph.  If
$\alpha(G)\leq m+k+1$, then $G$ has a spanning tree with at most $k$
branch vertices.
\end{corollary}

In 1998, Broersma and Tunistra gave the following degree sum condition for a graph to have a spanning tree with at most $l$ leaves.
\begin{theorem}[{\cite[Broersma and Tuinstra]{BT98}}]\label{t2}
	Let $G$ be a connected graph and let $l\geq2$ be an integer. If
	$\sigma_2(G)\geq |G|-l+1$, then $G$ has a spanning tree with at most
	$l$ leaves.
\end{theorem}
Motivating by Theorem \ref{t1}, a natural question is whether we can find sharp sufficient conditions of $\sigma_{l+1}(G)$ for a connected graph $G$ having a few leaves or branch vertices. This question is still open. But, in certain graph classes, the answers have been determined.

For a positive integer $r$, a graph is said to be \emph{$K_{1,r}$-free} if it does not contain $K_{1,r}$ as an induced subgraph. A $K_{1,3}$-free graph is also called a \emph{claw-free} graph. 

For the case of claw-free graphs, Gargano et al. proved the following.
\begin{theorem}[{\cite[Gargano et al.]{GHHSV}}]\label{thm1}
	Let $k$  be a non-negative integer and let $G$ be a connected claw-free graph of order $n$. If $\sigma_{k+3}(G)\geq n-k-2$, then $G$ has a spanning tree with at most $k$ branch vertices.
\end{theorem}
In 2020, Gould and Shull proved the following theorem which was a conjecture proposed by Matsuda et al. in \cite{MOY}. 
 \begin{theorem}[{\cite[Gould and Shull]{GS}}]\label{GS} Let $k$ be a non-negative interger and let $G$ be a connected claw-free graph of order $n.$ If $\sigma_{2k+3}(G)\geq n-2$, then $G$ has a spanning tree with at most $k$ branch vertices.
\end{theorem}
On the other hand, Kano et al. gave a sharp sufficient condition for a connected graph to have a spanning tree with few leaves.
\begin{theorem}[{\cite[Kano et al.]{KKMOSY12}}]\label{Kano12}
	Let $k$ be a non-negative integer and let $G$ be a connected claw-free graph of order $n$. If $\sigma_{k+3}(G)\geq n-k-2$, then $G$ has a spanning tree with at most $k+2$ leaves.
\end{theorem}
We note that the author \cite{H1} also introduced a new proof of Theorem \ref{Kano12} based on the techniques of Gould and Shull in \cite{GS}.

For connected $K_{1,4}$-free graphs, Kyaw~\cite{Ky09,Ky11} obtained
the following sharp results.

\begin{theorem} [{\cite[Kyaw]{Ky09}}]
	Let $G$ be a connected $K_{1,4}$-free graph with $n$ vertices. If
	$\sigma_4(G)\geq n-1$, then $G$ contains a spanning tree with at
	most $3$ leaves.
\end{theorem}

\begin{theorem} [{\cite[Kyaw]{Ky11}}] \label{Ky11}
	Let $G$ be a connected $K_{1,4}$-free graph with $n$ vertices.
	\begin{itemize}
		\item [$($i$)$]
		If $\sigma_3(G)\geq n$, then $G$ has a hamiltonian path.
		\item [$($ii$)$]
		If $\sigma_{m+1}(G)\geq n-\frac{m}{2}$ for some integer $m\geq 3$,
		then $G$ has a spanning tree with at most $m$ leaves.
	\end{itemize}
\end{theorem}
Regarding the existence of a spanning tree with a bounded number of branched vertices in a connected graph, Flandrin et al. proposed the following conjecture.

\begin{conjecture}[{\cite[Flandrin et al.]{FKKLR}}]\label{conj1}
	Let $k$ be a positive interger and let  $G$ be a connected  graph of order $n$. If $\sigma_{k+3}(G)\geq n-k$, then $G$ has a spanning tree with at most $k$ branch vertices. 
\end{conjecture}
Recently, Hanh gave a proof for Conjecture \ref{conj1} in the case graphs are $K_{1,4}$-free.
\begin{theorem}[{\cite[Hanh]{Hanh}}]\label{Hanh} Let $k$ be a positive interger and let
	$G$ be a connected $K_{1,4}$-free graph of order $n.$ If $\sigma_{k+3}(G)\geq n - k,$ then
	$G$ has a spanning tree with at most $k$ branch vertices.
\end{theorem}

For the $K_{1,5}$-free graphs, some results were obtained as follows.
\begin{theorem} [{\cite[Chen et al.]{CHH}}]
	Let $G$ be a connected $K_{1,5}$-free graph with $n$ vertices. If
	$\sigma_5(G)\geq n-1$, then $G$ contains a spanning tree with at
	most $4$ leaves.
\end{theorem}
\begin{theorem} [{\cite[Hu and Sun]{HS}}] Let $G$ be a connected $K_{1,5}$-free graph with $n$ vertices.
	If $\sigma_6(G) \geq  n - 1,$ then $G$ contains a spanning tree with at most $5$ leaves.
\end{theorem}

Moreover, many researchers have also studied the degree sum conditions for graphs to have
spanning trees with a bounded number of branch vertices and leaves.
\begin{theorem}[{\cite[Nikoghosyan]{N}, \cite[Saito and Sano]{SS}}]\label{t3}  Let $k\geq 2$ be an integer. If a connected graph $G$ satisfies $\deg_G(x) + \deg_G(y) \geq |G|-k+1$ for every two non-adjacent vertices $x, y \in V (G)$, then $G$ has
	a spanning tree $T$ with $|L(T)| + |B(T)| \leq k + 1.$
\end{theorem}
In 2019, Maezawa et al. improved the previous result by proving the following
theorem.
\begin{theorem}[{\cite[Maezawa et al.]{MMM}}]  Let $k\geq 2$	be an integer. Suppose that a connected graph $G$ satisfies $\max\{\deg_G(x), \deg_G(y)\} \geq \dfrac{|G|-k+1}{2}$ for every two non-adjacent vertices $x, y \in V (G)$, then $G$ has
	a spanning tree $T$ with $|L(T)| + |B(T)| \leq k + 1.$
\end{theorem}

In this paper, we study the spanning tree with a bounded number of leaves and branch vertices for the case of $K_ {1,4}$-free graph. In particular, our main result is the following.

\begin{theorem}\label{thm-main}
	Let $k, m$ be two non-negative intergers {\rm ($ m \leq k+1$)} and let  $G$ be a connected $K_{1,4}$-free graph of order $n$. If $\sigma_{m+2}(G)\geq n-k$, then $G$ has a spanning tree with at most $m+k+2$ leaves and branch vertices. 
\end{theorem}
\section{Applications of the main result}
In this section, we introduce some applications of Theorem \ref{thm-main}.\\

When $m=0,$ we have the following corollary which is a particular case of Theorem \ref{t3} if graphs are $K_{1,4}$-free.
\begin{corollary}
	Let $k$ be a possitive interger and let  $G$ be a connected $K_{1,4}$-free graph of order $n$. If $\sigma_{2}(G)\geq n-k$, then $G$ has a spanning tree with at most $k+2$ leaves and branch vertices.
	\end{corollary}

When $m=k+1,$ we state the following result.
\begin{theorem}\label{thm-main1}
	Let $k$ be a non-negative interger and let  $G$ be a connected $K_{1,4}$-free graph of order $n$. If $\sigma_{k+3}(G)\geq n-k$, then $G$ has a spanning tree with at most $2k+3$ leaves and branch vertices.
\end{theorem}
We may show that Theorem \ref{Ky11} ($i$) and the following theorem as corollaries of Theorem \ref{thm-main1}.
\begin{theorem}[{\cite[Mom\`{e}ge]{M}}]\label{thm-main1a}
	Let $G$ be a connected $K_{1,4}$-free graph of order $n$. If $\sigma_{2}(G)\geq \dfrac{2}{3}n,$ then $G$ has a Hamiltonian path.
\end{theorem}
Indeed, it follows from the assumptions of Theorem \ref{thm-main1a} we obtain that $\sigma_{3}(G)\geq \dfrac{3}{2}\sigma_{2}(G) \geq n$ (that also satisfies the assumption of Theorem \ref{Ky11} ($i$)). Now, using Theorem \ref{thm-main1} with $k=0$ and $m=1$ we conclude that $G$ has a spanning tree $T$ with at most $3$ leaves and branch vertices. If $|L(T)|=3$ then $|B(T)|\geq 1,$ this is a contradiction. Then $|L(T)|\leq 2,$ this mean that $T$ is a path. Therefore, $G$ has a Hamiltonian path.\\

Moreover, we note that if the tree $T$ has at most $2k+3$ leaves and branch vertices then $T$ has at most $k$ branch vertices. So Theorem \ref{thm-main1} is an improvement of Theorem \ref{Hanh}. Then we give an affirmative answer for Conjecture \ref{conj1} in the case of $K_{1,4}$-free graphs with a new approach.\\

We end this section by constructing an example to show that
the conditions of Theorem \ref{thm-main1}
is sharp. Let $k, p$ be positive integers. Let $P=x_1x_2 ... x_{k+1}$ be a path. Let $D_0, $ $D_1, ..., $ $D_{k+1}, $ $ D_{k+2}$ be copies of the complete graph $K_p$ of order $p$.
For each $i\in\{1, 2, ..., k+1\},$ join $x_i$ to all vertices of the graph $D_{i}$, join $x_1$ to all vertices of the graph $D_0$ and join $x_{k+1}$ to all vertices of the graph $D_{k+2}$. Then the resulting graph $G$ is a $K_{1,4}-$free graph. On the other hand, we have $|G|=n=k+1+(k+2)p$ and $\sigma_{k+3}(G)= n-k-1$,
but $G$ has no spanning tree with at most $2k+3$ leaves and branch vertices.
\section{Definitions and Notations}
In this section, we recall some definitions which need for the proof of main results.
\begin{definition}[\cite{GS}]
	Let $T$ be a tree. For any two vertices of $T$, say $u$ and $v$, are joined by a unique path, denoted $P_T[u,v]$.  We also denote $\{u_v\}=V(P_T[u,v])\cap N_T(u)$ and $e_v$ as the vertex incident to $e$ in the direction toward $v$.
\end{definition}

\begin{definition}[\cite{GS}]
	Let $T$ be a spanning tree of a graph $G$ and let $v\in V(G)$ and $e\in E(T)$. Denote $g(e,  v)$ as the vertex incident to $e$ farthest away from $v$ in $T$. We say $v$ is an {\it oblique neighbor} of $e$ { with respect to $T$} if $vg(e,  v)\in E(G)$. Let $X\subseteq V(G)$. The edge $e$ has an {\it oblique neighbor in the set} $X$ if there exists a vertex of $X$ which is an  oblique neighbor of $e$  with respect to $T$.
\end{definition}
 
 \begin{definition}[\cite{GS}]
	Let $T$ be a spanning tree of a graph $G$. Two vertices are {\it pseudoadjacent} with respect to $T$ if there is some $e\in E(T)$ which has them both as oblique neighbors. Similarly, a vertex set is {\it pseudoindependent} with respect to $T$ if no two vertices in the set are pseudoadjacent with respect to $T$.
\end{definition}

\begin{definition}
	Let $T$ be a tree with $B(T) \neq \emptyset,$ for each a vertex $x \in L(T)$, set $y_x \in B(T)$ such that $\left( V(P_T[x,y_x])\setminus \{y_x\} \right) \cap B(T) = \emptyset$. We delete $V(P_T[x,y_x])\setminus \{y_x\}$ from $T$ for all $x \in L(T)$. The resulting graph is a subtree of $T$ and is denoted by $R\_Stem(T)$. It is also called the {\it reducible stem} of $T.$
\end{definition}
For two distinct vertices $v,w$ of $T$, we always define the
\emph{orientation} of $P_T[v,w]$ is from $v$ to $w$. If $x \in V(P_T[v,w])$, then $x^+$ and $x^-$ denote the successor and predecessor of $x$ on $P_T[v,w]$ if they exist, respectively. We refer to~\cite{Di05} for terminology and notation not defined here. 

\section{Proof of Theorem  \ref{thm-main}}


Suppose that $G$ has no spanning tree with at most total $k+m+2$ leaves and branch vertices. Choose some spanning $T$ of $G$ such that:\\
(C1) $|L(T)|$ is as small as possible.\\
(C2) $|R\_Stem(T)|$ is as large as possible, subject to (C1).\\

By the contrary hypotheses, we note that $|L(T)|+|B(T)|\geq k+m+3.$\\
If $|B(T)|=0,$ then $|L(T)|=2.$ So $|L(T)|+|B(T)|=2 < k+m+3.$ This is a contradiction. Hence, $|B(T)|\geq 1$ and, in particular, $B(T)\not=\emptyset.$\\
 On the other hand, we have
$$|L(T)| =2+\sum_{b\in B(T)}(\deg_T(b)-2)\geq 2+|B(T)|.$$
So 
 \begin{align*}
  &2|L(T)|\geq |L(T)|+2+|B(T)|\geq k+m+5\geq m-1+m+5=2m+4\\
 &\Rightarrow |L(T)|\geq m+2.
 \end{align*}

We now have the following claims.

\begin{claim}\label{claim11} 
	$L(T)$ is independent.
\end{claim}
\pf
	Assume that two leaves $s$ and $t$ are adjacent in $G$. Then $s$ has some nearest branch vertex $b$.
	Let $T'=T-\{bb_s\}+\{st\}.$ Then $T'$ is a spanning of $G$ satisfying $|L(T')|<|L(T)|,$ the reason is that either $T'$ has only one new leaf $b_s$ and $s, t$ are not leaves of $T'$ or $s$ is still a leaf of $T'$ but $T'$ has no new leaf and $t$ is not a leaf of $T'$. This contradicts to the condition (C1). So the claim holds. 
\qed

\begin{claim}\label{claim12} 
Let $b\in B(T)$ and $x\in N_T(b).$ For each vertex $s\in L(T),$ if $b\in V(P_T[s,x])$ then $sx\not\in E(G).$
\end{claim}
\pf
Assume that $sx\in E(G).$ Consider the spanning tree $T'=T-\{bx\}+\{sx\}.$ Hence, $|L(T')|<|L(T)|$ (since $s$ is not a leaf of $T'$), a contradiction with the condition (C1). So the claim is proved. 
\qed

\begin{claim}\label{claim13} 
	Let $b, r$ be two branch vertices of $T$ such that $V(P_T[b,r])\cap B(T) =\{b,r\}.$ Let $s$ be a leaf of $T.$ If $sx\in E(G)$ for some $x\in V(P_T[b,r])\setminus \{b\}$ then  $sx^-\not\in E(G).$ 
\end{claim}
\pf
Assume that there exists a vertex $x\in V(P_T[b,r])\setminus \{b\}$ such that $sx, sx^- \in E(G)$ (note that possibly $x^- =b$). Let $c$ be the nearest branch vertex of $s$. Consider the spanning tree $T'=T-\{xx^-,ss_c\}+\{sx,sx^-\}.$ If $s_c=c$ then $s$ is not a leaf of $T'.$ Hence, $|L(T')|<|L(T)|$, a contradiction with the condition (C1). Otherwise, $L(T')=L(T)$ and $|R\_Stem(T')|>|R\_Stem(T)|$ (since $s\in V(R\_Stem(T'))$), a contradiction with the condition (C2). This completes the proof of claim. 
\qed

\begin{claim}\label{claim14} 
Let $b, r$ be two branch vertices of $T$ such that $V(P_T[b,r])\cap B(T) =\{b,r\}.$ If $x\in V(P_T[b,r])\setminus \{b,r\}$ then  $|N(L(T))\cap \{x\}|\leq 1.$ 
\end{claim}
\pf
Assume that there exists a vertex $x\in V(P_T[b,r])\setminus \{b,r\}$ such that  $|N(L(T))\cap \{x\}|\geq 2.$ Then there are two vertices $s, t\in L(T)$ such that $xs, xt \in E(G).$ Without loss of generality, we may assume that $b \in V(P_{T}[s,x]).$ By Claim \ref{claim12}, we obtain $x^-\not=b.$ Since Claim \ref{claim11} and Claim \ref{claim13} hold, we have $st, sx^-,sx^+,tx^-,tx^+ \not\in E(G)$ (here $x^+$ can be $r$). Moreover, $G[x,x^-,x^+,s,t]$ is not $K_{1,4}$-free. Hence, we obtain $x^-x^+ \in E(G).$ Let $c$ be the nearest branch vertices of  $s.$ Consider the spanning tree $T'=T-\{xx^-,xx^+,cc_s\}+\{sx,tx,x^-x^+\}.$ Hence, $|L(T')|<|L(T)|,$ the reason is that either $T'$ has only one new leaf $c_s$ and $s, t$ are not leaves of $T'$ or $s$ is still a leaf of $T'$ but $T'$ has no new leaf and $t$ is not a leaf of $T'$. This contradicts to the condition (C1). \\
Therefore, Claim \ref{claim14} is proved. 
\qed

\begin{claim}\label{claim15} 
	$L(T)$ is pseudoindependent with respect to  $T$.
\end{claim}
\pf
	Suppose two leaves $s$ and $t$ are pseudoadjacent with respect to $T$. Then there exists some edge $e\in E(T)$ such that  $sg(e,s), tg(e,t)\in E(G)$. Let $b$ and $u$ be the nearest branch vertices of  $s$ and $t$, respectively. Consider two cases as follows:
	
	Case 1. Suppose $g(e,s)\not= g(e,t)$. Then $e_s=g(e,t)$ and $e_t=g(e,s)$, so $se_t, te_s\in E(G)$. Then $T'=T-\{e, bb_s\}+\{se_t, te_s\}$ violates (C1) since $T'$ has only one new leaf $b_s$ and $s, t$ are not leaves of $T'$ or $s$ is still a leaf of $T'$ but $T'$ has no new leaf and $t$ is not a leaf of $T'.$ So the case 1 does not happen.
	
	Case 2: Suppose $g(e,s)=g(e,t)$. Define $x:=g(e,s)=g(e,t)$. Then $e_s=e_t$ and denoted by vertex $z$. We have $xs, xt\in E(G)$. Since $s, t\in L(T)$ and $L(T)$ is independent, we have $x\notin L(T)$. Then there exists some vertex $y\in N_T(x)\setminus\{z\}.$
	
	If $sz\in E(G)$ then we consider the spanning tree $T'=T-\{bb_s,  e\} + \{sz,  tx\}.$ It follows from Claim \ref{claim12} that $z \not\in B(T).$ Hence $|L(T')|< |L(T)|$ (since two leaves  $s$ and $t$ are lost while $b_s$ is gained or $s$ is still a leaf of $T'$ but $T'$ has no new leaf and $t$ is not a leaf of $T'$). So $sz\notin E(G)$. The same argument gives $tz\notin E(G)$.
	
	If $sy\in E(T)$ then the spanning tree $T'=T-\{uu_t,  e\} + \{sy,  tx\}$ violates (C1) (since two leaves  $s$ and $t$ are lost while $u_t$ is gained or $t$ is still a leaf of $T'$ but $T'$ has no new leaf and $s$ is not a leaf of $T'$). So $sy\notin E(G)$. The same argument gives $ty\notin E(G)$.
	
Now, since $G[x, y, z, s, t]$ is not $K_{1,4}$-free and  $st, sz, sy, tz, ty\notin E(G)$, we obtain $yz\in E(G).$ Then the spanning tree $T'=T-\{e, xy, bb_s\}+\{sx,tx,yz\}$ violates (C1), the reason is that either $T'$ has only one new leaf $b_s$ and $s, t$ are not leaves of $T'$ or $s$ is still a leaf of $T'$ but $T'$ has no new leaf and $t$ is not a leaf of $T'$.\\
 The claim \ref{claim15} has been proven.
\qed

\begin{claim}\label{claim16}  
	For each pair branch vertices $ b, r \in B(T) $ such that $V(P_T[b,r])\cap B(T)=\{b,r\}$, there exists some edge $e\in E(P_T[b,r])$ which has no oblique neighbor in the set $L(T)$.
\end{claim}
\pf
We consider three cases as follows.

Case 1. $V(P_T[b,r])=\{b,r\}.$ By Claim \ref{claim12} we choose $e=br$. 

Case 2. $V(P_T[b,r])\not=\{b,r\}.$ On $P_T[b,r]$ we set $x=b^+\not= r.$ Assume that there doesn't exist edge in $E(P_T[b,r])$ which has no oblique neighbor in the set $L(T).$ Hence both of $e=bx, f=xx^+$ (note that possibly $x^+=r$) have oblique neighbors in $L(T).$ Then there exist $s, t \in L(T)$ such that $sg(f,s), tg(e,t)\in E(G).$ 

By Claim \ref{claim12} we obtain that $g(e,t)=b.$ If $g(f,s)=x$ then $s\not=t$ (by Claim \ref{claim13}). Let $c$ be the nearest branch vertices of  $s.$ Consider the spanning tree $T':= T-\{e,cc_s\}+\{tb,sx\}.$ Hence, $|L(T')|<|L(T)|,$ the reason is that either $T'$ has only one new leaf $c_s$ and $s, t$ are not leaves of $T'$ or $s$ is still a leaf of $T'$ but $T'$ has no new leaf and $t$ is not a leaf of $T'$. This contradicts to the condition (C1). This implies $g(f,s)\not=x.$ Then, $g(f,s)=x^+.$\\
Since $b\in B(T)$, there exists some vertex $y\in N_T(b)\setminus\{x, b_s\}.$ By Claims \ref{claim12}-\ref{claim13}, we have $tb_s,ty,tx \not\in E(G).$ Combining with $G[b,x,b_s,y,t]$ is not $K_{1,4}$-free we obtain either $xy\in E(G)$ or $xb_s\in E(G)$ or $yb_s\in E(G).$

If $xy\in E(G)$ or $xb_s\in E(G)$ we consider the spanning tree 
$$T':=\left\{\begin{array}{ll}T-\{bx, by\}+\{bt,xy\},
& \;\mbox{ if } xy\in E(G),\\
T-\{bx, bb_s\}+\{bt,xb_s\}, & \;\mbox{ if } xb_s\in E(G).
\end{array}\right.$$
Then $|L(T')|< |L(T)|$ ($t$ is not a leaf of $T'$). This contradicts to the condition (C1).

 If $yb_s\in E(G)$ then the spanning tree $T':=T-\{by, bb_s,xx^+\}+\{bt,sx^+,yb_s\}$ violates the condition (C1), the reason is that $T'$ has only one new leaf $x$ and $s, t$ are not leaves of $T'.$
 
Therefore, Claim \ref{claim16} is proved.
\qed
\begin{claim}\label{claim17}  
	In the graph $G,$ there exists an independent set $S$ such that $|S|=m+2$ and there are at least $k$ distinct edges of $T$ which has no oblique neighbor in the set $S$.
\end{claim}
\pf
Since $|L(T)|\geq k+3,$ let $S$ be a subset in $L(X)$ such that $|S|=m+2.$ For each $x \in L(T)\setminus S,$ let $e$ be the edge of $T$ incident to $x.$ Then $x$ is an oblique neighbor of $e$ with respect to $T.$ Combining with Claim \ref{claim15} we obtain that $e$ has no oblique neighbor in the set $S$. Hence, there are at least $|L(T)|-m-2$ edges in $E(T)\setminus E(R\_Stem(T))$ which have no oblique neighbor in the set $S.$

On the other hand, consider the tree $H$ with vertex set $V(H)=B(T)$ and edge set $E(H)=\{br |\, b, r\in V(H)\, \text{and}\, V(P_T[b,r])\cap B(T)=\{b,r\}\}$ (here $E(H)$ can be an empty set if $|B(T)|=1$). By Claim \ref{claim16}, the number of edges of $R\_Stem(T)$ which has no oblique neighbor in the set $L(T)$ is greater than or equal to the number of edges of $H.$ Hence, there are at least $|E(H)|$ edges in $E(R\_Stem(T))$ which have no oblique neighbor in the set $S.$

Set $h$ to be the number of edges of $T$ which has no oblique neighbor in the set $S.$ By the arguments mentioned above, we conclude that
\begin{align*}
h&\geq |L(T)|-m-2 + |E(H)|=|L(T)|-m-2+|V(H)|-1 \\
&= |L(T)|-m-2+|B(T)|-1=|L(T)|+|B(T)|-m-3 \geq k.
\end{align*}

This completes the proof of Claim \ref{claim17}.
\qed

For any $v, x\in V(T)$, we have $vx\in E(G)$  if and only if  $v$ is an oblique neighbor of $xx_v$.  Therefore, the number of edges of $T$ with $v$ as an oblique neighbor equals the degree of $v$ in $G$. Combining with  Claim \ref {claim11}, Claim \ref{claim15} and Claim \ref {claim17}, we obtain that
\[\sigma_{k+3}(G) \leq \sum_{x\in S}\deg_G(x) \leq |E(T)|-k=|V(T)|-1-k=n-1-k,\]
which contradicts the assumption of Theorem \ref {thm-main}. The proof of Theorem \ref {thm-main} is completed. \qed


\end{document}